\documentclass[12pt]{article}
\usepackage{graphicx}
\usepackage{latexsym,amssymb,amsmath,amscd,amsfonts}
\topmargin - 2cm
\oddsidemargin - 0.1cm
\evensidemargin - 5cm
\newtheorem{theorem}{Theorem}[section]
\newtheorem{lemma}[theorem]{Lemma}
\newtheorem{cor}[theorem]{Corollary}

\newcommand{\Z}{\mathbb{Z}}
\newcommand{\vf}{\varphi}
\newcommand{\ve}{\varepsilon}

\newcommand{\supp}{\mbox{supp}}

\newcommand{\ds }{ \displaystyle }
\newcommand{\re}{\mathbb{R}}

\begin{document}

\title{\bf Continuous time random walk models for fractional space-time diffusion equations}

\author{Sabir Umarov}
\date{}
\maketitle

\begin{center}
{\it University of New Haven, Department of Mathematics, \\300 Boston Post Road, West Haven, CT 06516, USA}
\end{center}
\vspace{10pt}

\begin{abstract}
In this paper continuous time random walk models approximating fractional space-time diffusion processes are studied. Stochastic processes associated with the considered equations represent time-changed processes, where the time-change process is a L\'evy's stable subordinator with the stability index $\beta \in (0,1).$  In the parer the convergence of constructed CTRWs to time-changed processes associated with the corresponding fractional diffusion equations are proved using a new analytic method. 
\end{abstract}

%
\vspace{2pc}

%
%
%
%

\section{Introduction}\label{sec:1}

By definition, a continuous time random walk (CTRW) is a random walk subordinated to a renewal process. More precisely this means that CTRW comprises of two i.i.d. sequences of random variables (vectors), $X_n$ expressing jumps, and ${T}_n$ expressing waiting times between successive jumps. 
 CTRW model was introduced by Montrol and Weiss \cite{MontrolWeiss} in 1965. CTRW processes have broad applications in various fields (see, e.g. \cite{MetzlerKlafter,GorenfloScalas}) and have been extensively investigated recent years by many authors; see \cite{MetzlerKlafter}-\cite{Fa} 
and references therein. If the sequences $X_n$ and $T_n$ are independent, then STRW is called decoupled. We would like to single out two papers \cite{GorenfloMainardi, Liu}, contents of which are close to the present paper, and which discuss a continuous time random walk approximation of the stochastic process associated with the fractional order diffusion equation of the form
\begin{equation}
\label{eqn_010}
D_{\ast}^{\beta} u(t,x) = \frac{1}{2} \frac{\partial^2 u(t,x)}{\partial x^2}, \quad t>0, \, x \in \re.
\end{equation}
Here $D_{\ast}^{\beta}$ is the fractional derivative of order $\beta \in (0,1)$ in the sense of Caputo-Djrbashian:
\begin{equation}
\label{Cap}
D_{\ast}^{\beta} g(t)=\frac{1}{\Gamma(1-\beta)}\int_0^t \frac{g^{\prime}(\tau)}{(t-\tau)^{\beta}}d\tau.
\end{equation}
One can always assume that the process starts from the origin, which in terms of $u(t,x),$ is written in the form \begin{equation}
\label{in_cond}
u(0,x)=\delta_0(x),
\end{equation} 
where $\delta_0$ is the Dirac delta function concentrated at $0.$ These papers suggest two distinct discretizations of equation (\ref{eqn_010}). Both discretizations lead to the unified form
\begin{equation}
  D_{\ast}^{\beta} u^n_j \approx 
    a(\tau)                                              
  \left( 
     u^{n+1}_j - c_1 u^n_j - \sum_{m=2}^n c_m u^{n+1-m}_j - \gamma_n u^0_j 
  \right), \quad n=1,2,\dots,
  \label{eqn_011} 
\end{equation}
where $u_j^n=u(t_n,x_j),$ $t_n=n\tau, \, \tau>0, \, x_j = hj=(hj_1,\dots, h j_d),$ and 
$
u_j^0=1,
$
if $j=0=(0,\dots,0),$ and $u_j^0=0,$ if $j \ne 0.$ In \cite{GorenfloMainardi} the Gr\"unwald-Letnikov approximation is used for discretization of $D_{\ast}^{\beta} u(t,x),$ and parameters $a(\tau), \gamma_n,$ and  $c_k$ are defined as 
\begin{equation}
\label{c_1}
a(\tau) =a(\tau,\beta)=\frac{1}{\tau^{\beta}}, \\
c_k = c_k(\beta)= (-1)^{k} {\beta \choose k}, \, k=1,2,\dots,
\end{equation}
\begin{equation}
\label{c_3}
\gamma_n =\gamma_n(\beta)=\sum_{k=0}^n c_k, \, n=1,2,\dots.
\end{equation}
Paper \cite{Liu} uses a special quadrature for the discretization of $D_{\ast}^{\beta} u(t,x)$ with parameters
\begin{equation}
\label{c_4}
a(\tau) =\frac{1}{\tau^{\beta}\Gamma(2-\beta)}, 
\gamma_m  =  (m+1)^{(1-\beta)}-m^{(1-\beta)},  m=0, 1, \ldots,n, 
\end{equation}
\begin{equation}
\label{c_6}
  c_0  =1, \, c_k  = \gamma_{k-1} - \gamma_k,   k=1,\ldots,n,
\end{equation}
for an arbitrary fixed $n.$ For the right hand side of (\ref{eqn_010}) the second order symmetric finite difference approximation is used. These discretisations lead to an explicit scheme
$
u_j^{n+1} = F(u_j^{0},\dots,u_j^{n}),
$
where $F$ is a linear function of its arguments. 
The positiveness and conservativeness conditions, as well as 
the stability and convergence issues are studied in these papers using numerical methods. 

In our paper we develop an analytic method for convergence of the CTRW to the stochastic process associated with the distributed order diffusion equation 
\begin{equation*}
D_{\ast}^{\beta} u(t,x)=\int_{0}^2 \mathbb{D}_0^{\alpha} u(t,x) d\rho(\alpha), \quad t>0, \, x \in \re^d, \\
\end{equation*}
with initial condition (\ref{in_cond}). Here
$\mathbb{D}_0^{\alpha}, \, 0<\alpha <2,$ is a hypersingular integral operator defined in the form 
\begin{equation}
\label{equation_010}
\mathbb{D}_0^{\alpha}f(x)=  b(\alpha)\int\limits_{{\re}_y^n}\frac{\Delta_y^2
f(x)}{|y|^{n+\alpha}} dy,
\end{equation}
with $\Delta_y^2,$  the second order centered finite  difference
in the $y$ direction and $b(\alpha)$ is the normalizing constant
\[
b(\alpha)= \frac{\alpha \Gamma(\frac{\alpha}{2})
\Gamma(\frac{n+\alpha}{2}) \sin \frac{\alpha \pi}{2} } {2^{2-\alpha}
\pi^{1+n/2}},
\]
and $\rho$ is a finite measure with  support in the interval $[0,2].$ In this paper we will use only the Gr\"unwald-Letnikov approximation for $D_{\ast}^{\beta}.$ The method developed in this paper is a modification of the method used  in \cite{UmarovSteinberg,Andries} and provides a rigorous justification of the convergence of constructed CTRWs to their associated limit stochastic processes.


\section{Preliminaries and auxiliaries} 
\label{Preliminaries}

Let $\rho$ be a finite measure with the support  $\supp \, \rho \subseteq [0,2].$ Denote by
 $\mathbb{SS}$ the class of $d$-dimensional 
$(\mathcal{F}_t)$-semimartingales $Z_t, \, Z_0=0,$ whose characteristic function is
given by
\begin{equation} \label{st_mix_010}
 \mathbb{E} \bigl[e^{i \xi Z_t}\bigr]= \exp \Big\{-t \int\limits_0^2 |\xi|^{\alpha} d \rho(\alpha) \Big\}, ~ \xi \in \re^d.
\end{equation}
If $f_{Z_t}(x), \, x \in \re^d,$ is the density function of the process  $Z_t,$ then equation (\ref{st_mix_010}) can be expressed as the Fourier transform $\mathbb{E} \bigl[e^{i \xi Z_t}\bigr] = F[f_{Z_t}](\xi).$ The class 
$\mathbb{SS}$
obviously contains L\'evy's symmetric $\alpha$-stable processes and all mixtures of their finitely
many independent representatives.  For the process $Z_t \in \mathbb{SS}$
corresponding to a finite measure $\rho$, we use the notation
$Z_t=X_t^{\rho}$ to indicate this correspondence.

Suppose $X^{\rho}_t \in \re^d$ is a stochastic process obtained by mixing of independent L\'evy's $S\alpha S$-processes with a mixing measure $\rho,\, \supp \, \rho \subset (0,2).$
Then its associated FPK equation has the form
\begin{equation} 
\label{equation} 
\frac{\partial u(t,x)}{\partial t}= \int\limits_0^2 \mathbb{D}_0^{\alpha}u(t,x)d\rho(\alpha), \quad t>0, \, x \in \re^d,
\end{equation}
where $\mathbb{D}_0^{\alpha}$ is defined in (\ref{equation_010}). Indeed, it is seen from (\ref{st_mix_010}) that the L\'evy symbol of 
$X_t^{\rho}$ equals 
\begin{equation}
\label{fr_psi_010}
\Psi(\xi)=-\int_0^2 |\xi|^{\alpha} d \rho.
\end{equation}
On the other hand, due to the known relationship $F[D_0^{\alpha} \vf](\xi)=-|\xi|^2F[\vf](\xi)$ (see, e.g. \cite{SKM}), the Fourier transform of  the right hand side of (\ref{equation}) is
\[
F\Big[ \int_0^2 D_0^{\alpha} \vf (x) d\rho(\alpha)\Big] =\int_0^2 F[D_0^{\alpha} \vf ](\xi) d \rho(\alpha) 
\]
\[
= \Big(-\int_0^2 |\xi|^{\alpha} d \rho \Big)\, F[\vf](\xi) = \Psi(\xi) F[\vf](\xi), \quad \xi \in \re^d.
\]
Hence, the symbol of the pseudo-differential operator on the right hand side of (\ref{equation}) coincides with $\Psi(\xi).$ This means that the FPK equation associated with $X_t^{\rho}$ is given by equation (\ref{equation}). Since $X^{\rho}_0=0$ the function $u(t,x)$ in (\ref{equation}) satisfies the initial condition $u(0,x)=\delta_0(x).$
Using the Fourier transform technique, one can verify that the solution of (\ref{equation}) satisfying the latter initial condition can be expresses in the form
\begin{equation}
\label{fr_sln_090}
G^{\rho}(t,x)=\frac{1}{(2\pi)^n} \int_{\re^n} e^{t \Psi(\xi)-i x \xi} d\xi.
\end{equation}

We note also that a pseudo-differential operator $A(D)$ and its symbol $\sigma_{A}(\xi)$ are connected through the relation \cite{Hormander} $\sigma_A(\xi)=e^{-ix \xi} A(D) e^{ix\xi}$ valid for all $x.$ In particular, if $x=0:$
\begin{equation}
\label{symbol}
\sigma_A(\xi)= \big( A(D) e^{ix\xi} \big)_{|x=0}, \quad \xi \in \re^d.
\end{equation}

Let $X^{\rho}_t \in \mathbb{SS}$ and $Y_t=X^{\rho}_{W_t}$ be the time-changed process, where $W_t$ is the inverse to a
$\beta$-stable subordinator (see \cite{Sato}). Then the density of $Y_t$  solves the fractional FPK type equation \cite{HKU}
\begin{equation}
\label{ctrw_010}
D_{\ast}^{\beta} u(t,x)= \int_0^2 \mathbb{D}_0^{\alpha} u(t,x) d\rho, \quad t>0, \, x \in \re^d,
\end{equation}
with a finite measure $\rho, \, \supp \, \rho \subset (0,2],$ and initial condition (\ref{in_cond}).  The unique solution of equation (\ref{ctrw_010}) satisfying initial condition (\ref{in_cond}) is given by \cite{UmarovGorenflo}
\[
u(t,x)= E_{\beta} (\Psi(D)t^{\beta} ) \delta_0(x) = \frac{1}{(2\pi)^d} \int_{\re^d} E_{\beta}[\Psi(\xi)t^{\beta}]  e^{-i x \xi} d\xi,
\]
where $E_{\beta}(z)$ is the Mittag-Leffler function \cite{ML}.
Thus, the Fourier transform of $u(t,x)$ in the variable $x$ is $F[u(t,\cdot)](\xi)=E_{\beta}(\Psi(\xi)t^{\beta}).$
Moreover, it follows from the propertyies of the Mittag-Leffler function that the Fourier-Laplace transform of $u(t,x)$ is
\[
L[F[u](\xi)](s)=\frac{s^{\beta-1}}{s^{\beta}+(-\Psi(\xi))}, \quad s>0, \, \xi \in \re^d.
\]
We will use the Gr\"unwald-Letnikov approximation for discretization of the Caputo-Djrbashian fractional derivative in equation (\ref{ctrw_010}). By definition, the backward Gr\"unwald-Letnikov fractional derivative of order $\beta$ of a function, defined on an interval $(a,b),$ is
\begin{equation}
\label{GL}
_a \mathcal{D}_t^{\beta} f(t)=\lim_{h \to 0} \frac{1}{\tau^{\beta}} \sum_{m=0}^{\lfloor \frac{t-a}{\tau} \rfloor} (-1)^m {\beta \choose m} f(t-m\tau).
\end{equation}
The Gr\"unwald-Letnikov fractional derivative (\ref{GL}) can be used for approximation of the Caputo-Djrbashian derivative (\ref{Cap}). Namely, taking $\tau =(t-a)/n$ and $t_k=a+k\tau, k=0,\dots,n,$ in (\ref{GL}), one has  the approximation \cite{Gorenflo}
\begin{equation}
\label{Cap=GL}
D_{\ast}^{\beta} f(t_n) \approx \frac{1}{\tau^{\beta}} \sum_{m=0}^{n}(-1)^m {\beta \choose m} \Big(f(t_{n-m})-f(a)\Big).
\end{equation}
with the order of accuracy $O(\tau).$
Moreover, if $a =-\infty$, then the finite sum in (\ref{GL})  becomes an infinite series, that is
\begin{equation}
\label{GL_1}
_{\infty} \mathcal{D}_t^{\beta} f(t)=\lim_{\tau \to 0} \frac{1}{\tau^{\beta}} \sum_{m=0}^{\infty} (-1)^m {\beta \choose m} f(t-m\tau),
\end{equation}
convergent for functions satisfying the asymptotic behavior 
\begin{equation}
\label{cond}
f(t) =O(|t|^{-(1+\beta+\ve)}), \quad t \to -\infty,
\end{equation} 
for some $ \ve>0.$ It is also known \cite{SKM} that in the class of such functions $_{\infty} \mathcal{D}_t^{\beta} f(t)$ coincides with the Liuoville-Weyl backward fractional derivative of order $\beta, $ defined as 
\[
_{-\infty} \mathfrak{D}_t^{\beta} f(t)= \frac{1}{\Gamma(1-\beta) } \frac{d}{dt}\int_{-\infty}^t \frac{f(\tau)}{(t-\tau)} d\tau, \quad (0<\beta<1).
\]
For example, the function $f(t)=e^{st}, \, s>0,$ defined on $(-\infty, b], \, 0<b<\infty,$ obviously satisfies condition (\ref{cond}). It can be readely veryfied that  $_{-\infty} \mathfrak{D}_t^{\beta} e^{st}=s^{\beta} e^{st}.$
Hence, taking $t=0,$ one has
\begin{equation}
\label{ex}
\left( _{\infty} \mathcal{D}_t^{\beta} e^{st} \right)_{|t=0}= \left(_{-\infty} \mathfrak{D}_t^{\beta} e^{st}\right)_{|t=0} = s^{\beta}.
\end{equation}
This example will be used in in the proof of the main result.

\begin{lemma} \label{lem}
Let   $0<\alpha<2$ and
\begin{equation} \label{pmf_p_k}
p_{k} =
\left\{
             \begin{array}{ll}
             {\ds b(\alpha) |k|^{-(d + \alpha)}}, & k \neq 0, \\ & \\
             {\ds -b(\alpha) \sum_{m \in {\Z}^d \backslash \{0\}}} \frac{1}{|m|^{d + \alpha}}, & k = 0,
             \end{array}
           \right.
\end{equation}
where $\Z^d$ is the $d$-dimensional integer lattice.
Then the characteristic function 
\[
\hat{p}_h(\xi)=\sum_{k \in \Z^d} p_k e^{ikh \xi},
\]
of the sequence $\{p_k\}_{k \in \Z^d},$ converges to $-|\xi|^{\alpha}/2$ as $h \to 0.$ 
\end{lemma}

{\it Proof.}  Let
$\hat{p}_h(\xi)$ be the characteristic function of $p_k$ defined in (\ref{pmf_p_k}), that is
\[
\hat{p}_h(\xi) = \sum_{k \in {\Z}^d} p_{k} e^{i hk \xi} 
\]
\[
= -b(\alpha) \sum_{m \in {\bf Z}^d \backslash \{0\}} \frac{1}{|m|^{d + \alpha}} + b(\alpha) \sum_{m \in {\bf Z}^d \backslash \{0\}} \frac{1}{|m|^{d + \alpha}} e^{imh\xi}
\]
\[
= -b(\alpha) \sum_{m \in {\bf Z}^d \backslash \{0\}} \frac{1}{|m|^{d + \alpha}} (1-e^{im h\xi} ).
\]
Further, one can easily verify that
\begin{equation*} \label{2-3_line}
\sum_{0 \ne k \in
{\Z}^d} \frac{1-e^{ik \xi h}}{|k|^{d+\alpha}}  = \sum_{0 \ne k \in
{\Z}^d} \frac{1-e^{-ik \xi h}}{|k|^{d+\alpha}}.
\end{equation*}
Due to the definition of the symmetric second finite difference of the function $e^{i x \xi}$  at the origin, this implies
\[
\hat{p}_h(\xi)=- b(\alpha)
\sum_{0 \ne k \in
{\Z}^d} \frac{1-e^{ik \xi h}}{|k|^{d+\alpha}}   =-b(\alpha){1 \over 2} \sum_{0 \ne k \in
{\Z}^d} \frac{2-e^{ik\xi h}+e^{-ik \xi h}}{|k|^{d+\alpha}} 
\]
\begin{equation*} 
=b(\alpha){1 \over 2} \sum_{0 \ne k \in
{\Z}^d} \frac{ (\Delta_{kh}^2 e^{i x \xi})_{|x=0} } {|k|^{d+\alpha}}.
\label{2-3_line_1}
\end{equation*}
Now letting $h \to 0,$ due to definition (\ref{equation_010}) of $\mathbb{D}_0^{\alpha}$ and relation (\ref{symbol}), we have
\[
\lim_{h\to 0}\hat{p}_h(\xi)={1 \over 2}b (\alpha)\int_{{\re}^d}
                     \frac{(\Delta_y^2 e^{ix\xi})_{x=0} }{|y|^{d+\alpha}} dy   = {1 \over 2}(D_0^{\alpha} e^{ix\xi})_{|_{x=0}} 
                     = - {1 \over 2} |\xi|^{\alpha},
\] 
as desired.

\begin{cor}
\label{cor}
Let   
\begin{equation} 
d_k := \left\{
             \begin{array}{ll}
             {\ds \frac{Q_k(h)}{|k|^d}}, & k \neq 0, \\ & \\
             {\ds - \sum_{k \neq 0} \frac{Q_k(h)}{|k|^d}}, & k = 0,
             \end{array}
           \right.
  \label{ctrw_eq:space_frac_der_0}
\end{equation}
where
\begin{equation}
\label{fr_rw_051}
Q_k(h)= 2\int_0^2 \frac{b(\alpha)d\rho(\alpha)}{h^{\alpha}|k|^{\alpha}}, \quad k \neq 0.
\end{equation}
Then the characteristic function  ${\hat{d}_h(\xi)}$ of the sequence $\{d_k\}_{k \in \Z^d},$ converges to $\Psi(\xi)$ as $h \to 0.$
\end{cor}

{\it Proof} \
We have
  \[
  \hat{d}(h \xi) =- \sum_{k\ne 0} \frac{Q_k(h)}{|k|^d}+ \sum_{k\ne 0} \frac{Q_k(h)}{|k|^d}e^{ik h \xi}
  \]
\[
  =  \sum_{k \ne 0} \frac{Q_k(h)}{|k|^d} (e^{ikh \xi}-1)
  ={1 \over 2} \sum_{k \ne 0} \frac{Q_k(h)}{|k|^d} \left(e^{ik h \xi} -2 +e^{ikh \xi} \right) 
  \]
\[
  = {1 \over 2} \sum_{k \ne 0} \frac{Q_k(h)}{|k|^d} \left(\Delta_{kh}^2 e^{i x \xi}  \right)_{|x=0}  
  = \int_0^2 b(\alpha) \left( \sum_{k \ne 0} \frac{\Delta_{kh}^2 e^{i x \xi}}{|kh|^{d+2}} h^d \right)_{|x=0} d\rho (\alpha) 
\]
\[
  =\int_0^2 \hat{p}_h(\xi) d\rho(\alpha),
  \]
where $\hat{p}_h(\xi)$ is the charactersitic function of the sequence $p_k$ in equation (\ref{pmf_p_k}).
  Letting $h \to 0,$ due to Lebesgue's dominated convergence theorem and Lemma \ref{lem}, 
 we obtain
  \[
  \lim_{h \to 0} \hat{d}(h\xi) 
  = -\int_0^2 |\xi|^2 d\rho(\alpha) = \Psi(\xi).
  \]

\section{Main results: CTRWs and their limit processes}

As is known \cite{HKU,HU,HKU2}, driving processes of stochastic differential equation associated with time-fractional Fokker-Planck-Kolmogorov (FPK)
equations appear to be time-changes of basic processes like Brownian
motion, L\'evy process, fractional Brownian motion, etc. Donsker's
theorem states that Brownian motion is the limit in the weak
topology of a scaled sum of a sequence of independent and
identically distributed (i.i.d.) random variables $\{X_j\},$ with
$X_1 \in L^2(P)$. This fact is important from the approximation
point of view since an approximation of the basic driving process
$B_t$ yields, under some conditions, an approximation of other
processes $X_t$ driven by $B_t.$ Natural approximants of
time-changed processes $B_W, \, L_W,$ etc., where $W$ is the inverse
to a stable subordinator, are CTRWs. 
A decoupled CTRW is defined by two independent sequences of 
random variables/vectors: one representing the sizes of jumps, the
other representing waiting times between successive jumps. More precisely, let 
$$
Y_1, Y_2, \dots, Y_n, \dots, \quad (Y_i \in \re^d),
$$ 
be a sequence of i.i.d. random vectors, and 
let 
$$
T_1, T_2, \dots, T_n, \dots, \quad (T_i \in \re_+)
$$ be an i.i.d.
sequence of positive real-valued random variables. 
Then 
$$S_n=Y_1+\cdots+Y_n$$ 
is the position after $n$ jumps, and
$$t_n=T_1+\cdots+T_n$$ 
is the time of the $n$th jump. Assume
that $S_0=0$ and $t_0=0$. The stochastic process 
\[
X_t= 
S_{N_t}= \sum_{i=1}^{N_t}Y_i,
\]
{where } 
$
N_t =\max\{n\ge0: t_n \le t\},
$
is called a \emph{continuous time random walk}.

Since $Y_t=X^{\rho}_{W_t}$ is a non-Markovian process, its approximating random walk also can
not have independent members. Therefore, transition probabilities split into
two different sets of probabilities:
\begin{enumerate}
\item \textit{non-Markovian transition probabilities,} which express a long
non-Markovian memory of past; and
\item \textit{Markovian transition probabilities,} which express transition
from positions at the previous time instant.
\end{enumerate}
%
Suppose that non-Markovian transition probabilities are given by
\[
c_{\ell}=(-1)^{\ell+1} {\beta \choose \ell}=\Big|{\beta \choose
\ell}\Big|, \,
\ell=1,\dots,n, 
\]
\begin{equation}
\gamma_n=\sum_{\ell=0}^n(-1)^{\ell}{\beta \choose \ell}, \label{cbn}
\end{equation}
and Markovian transition probabilities $\{p_k\}_{k\in \mathbb{Z}^n}$
are given by
\begin{equation}\label{trpr11}
p_{k} = \left\{ \begin{array}{ll}
          c_1-  \tau^{\beta} Q(h), & \mbox{if $k=0$;}
\\         \tau^{\beta} \frac{Q_k(h)}{|k|^{d}},  & \mbox{if $k \neq 0$,} \end{array} \right.
\end{equation}
where $Q_k(h), \, k \neq 0,$ is defined in (\ref{fr_rw_051}), and $Q(h)=\sum_{k \neq 0} Q_k(h)|k|^{-n}.$
Then CTRW, approximating the stochastic process $Y_t,$ can be interpreted in the following sense: the probability $q_j^{n+1}$ of the walker being at a site $x_j=jh, \, j \in \Z^d,$ at a
time $t_{n+1}$ is
\begin{equation}
\label{ctrw_000}
q_j^{n+1}=\gamma_n
q_j^0+\sum_{\ell=1}^{n-1}c_{n-\ell+1}q_j^{\ell}+\Big(c_1-\tau^{\beta}Q_0(h)\Big)q_j^n+\sum_{k
\ne 0}p_kq_{j-k}^n.
\end{equation}

To construct CTRW (\ref{ctrw_000}) one needs to discretize (\ref{ctrw_010}). Namely,
for the Caputo fractional derivative on the left-hand-side of (\ref{ctrw_010}) 
 we use the backward Gr\"unwald-Letnikov  discretization (\ref{Cap=GL}) in the form:
\begin{equation}
   D_{\ast}^{\beta} u^{n}_j = _0 \mathcal{D}_t^{\beta} u^{n}_j
   \approx   
  \sum_{m=0}^{n+1}  (-1)^{m}
  {\beta \choose m} \frac{u_j^{n+1-m}-u_j^0}{{\tau^{\beta}}}
   \label{ctrw_eq:time_frac_der_discr_0}
\end{equation}
where $u_j^n=u(t_n,x_j), \, n=0,1,\dots, \, j\in \Z^d,$ $x_j\in h\Z^d,$ and $t_n=n\tau, \, \tau >0.$ 
Using notations (\ref{cbn})
and rearranging terms, equation 
(\ref{ctrw_eq:time_frac_der_discr_0}) can be expressed in the form
\begin{equation}
  D_{\ast}^{\beta} u^{n}_j \approx 
  \frac{1}{ \tau^{\beta}} 
  \left( 
     u^{n+1}_j - c_1 u^n_j - \sum_{m=2}^{n+1} c_m u^{n+1-m}_j - \gamma_n u^0_j 
  \right).
  \label{ctrw_eq:time_frac_der_discr_caputo_0} 
\end{equation}

For  the right hand side of (\ref{ctrw_010}) we use the discretization 
\begin{equation}
  \Psi(D_x) u^n_j  \approx
  \sum_{k \in \Z^d} d_k u^n_{j-k},  
  d_k := \left\{
             \begin{array}{ll}
             {\ds \frac{Q_k(h)}{|k|^d}}, & k \neq 0, \\ & \\
             {\ds - \sum_{k \neq 0} \frac{Q_k(h)}{|k|^d}}, & k = 0,
             \end{array}
           \right.
  \label{ctrw_eq:space_frac_der_0}
\end{equation}
where $Q_k(h)$ is defined in (\ref{fr_rw_051}).
Setting the discretizations for the time and space-fractional 
derivatives in (\ref{ctrw_eq:time_frac_der_discr_caputo_0}) and (\ref{ctrw_eq:space_frac_der_0}) 
equal to each other, we get  
\begin{equation} \label{eq}
  \frac{1}{\tau^{\beta}} 
  \left( u^{n+1}_j - c_1 u^n_j - \sum_{m=2}^n c_m u^{n+1-m}_j 
         - \gamma_n u^0_j
  \right) = 
  \sum_{k \in \Z^d} d_k u^n_{j-k}.    
\end{equation}
Rearranging terms and solving for $u^{n+1}_j$ in equation (\ref{eq}), the following 
recursion equation is constructed, reconstituting CTRW (\ref{ctrw_000}):
\begin{eqnarray} 
  u^{n+1}_{j} 
  & = &  
  {\ds \gamma_n u_j^0 + \sum_{m=2}^{n} c_m u_j^{n+1-m} +
                   \sum_{k \in \Z^d} q_k u^n_{j-k}}, 
  \label{fr_eq:time_space_frac_discr_eqn_0} \\
  & & \nonumber \\ 
  q_k
  & = & 
  \left\{ \begin{array}{ll}
          {\ds \tau^{\beta} d_k =
           2\int_{0}^2
           \left( \frac{\tau^{\beta}}{h^{\alpha}} \right)
          \frac{d \rho(\alpha)}{|k|^{d+\alpha}}}, & k \neq 0 \\ 
          & \\
          {\ds c_1 - \sum_{k \neq 0} q_k}, & k = 0.
          \end{array} \right. \nonumber
\end{eqnarray}    
By construction, $u_j^0=1$ if $j=0=(0,\dots, 0),$ and $u_j^0=0$ otherwise.

The update $u^{n+1}_j$ in equation (\ref{fr_eq:time_space_frac_discr_eqn_0}) 
is determined by Markovian contributions
(those values of $u$ at time $t=t_n$) and non-Markovian 
contributions (those values of $u$ at times 
$t=\{t_0,t_1,\ldots,t_{n-1}\}$). 
The order of the time fractional derivative $\beta$
determines the effect that the non-Markovian transition
probabilities ($\gamma_n$ and $c_2,\ldots,c_m$)
has on $u^{n+1}_j$.  This effect can be measured by 
sum of all of the transition probabilities in equation
(\ref{fr_eq:time_space_frac_discr_eqn_0}):
\[
  {\ds \left( \gamma_n + \sum_{m=2}^n c_m \right)  + 
       \sum_{k \in \Z^d} q_{k} = 1}. 
  \label{eq:transprob_sum}
\]
where
\[
  \sum_{k \in \Z^d} q_{k} = (c_1 - q_0 ) + \sum_{k \neq 0} q_k =  c_1
  \hspace*{.75cm} \mbox{and} \hspace*{.75cm}
  \gamma_n + \sum_{m=2}^n c_m =  1 - c_1.
\]
As a result, when $\beta=1$ one has $c_1 = 1$,
$c_2 = \cdots = c_n = \gamma_n = 0,$ 
and hence, equation (\ref{fr_eq:time_space_frac_discr_eqn_0})
simply reduces to 
\begin{equation*}
\label{fr_eq_recursion}
u_j^{n+1}=\sum_{k \in {\Z}^d}p_k u_{j-k}^n, \,  j \in {\Z}^d, \, 
n =0,1, \dots.
\end{equation*}
with the transition probabilities 
\begin{equation}
\label{trpr111}
p_{k} = \left\{ \begin{array}{ll}
          1-  \tau 
          \sum_{m \neq 0} \frac{Q_m(h)}{|m|^{d}}, & \mbox{if $k=0$;}
\\  \tau \frac{Q_k(h)}{|k|^{d}}.
      & \mbox{if $k \neq 0$,} \end{array} \right.
\end{equation}

\begin{theorem}  \label{thm}
Let $0<\beta \le 1.$ Fix $t>0$ and let $h>0, \, \tau =t/n.$
Let $Y_j \in \mathbb{Z}^d, \, j \ge 1,$ be identically distributed
random vectors with the non-Markovian and Markovian transition
probabilities defined in (\ref{cbn}) and in (\ref{trpr11}),
respectively. Assume that
\begin{equation} \label{stabcond2}
\tau \le \Big(\frac{\beta}{Q(h)}\Big)^{{1 \over \beta}}.
\end{equation}
Then the sequence of random vectors ${S}_{n}=h{Y}_{1}+ ... + h {
Y}_{n},$ converges as
 $n \rightarrow \infty$ in law to
$X_t=Y_{W_t}$ whose probability density function is the solution to
equation (\ref{ctrw_010}) with the initial condition
$u(0,x)=\delta_0 (x).$
\end{theorem}

{\it Proof.}  
   Let $\hat{u}^n(\xi)$ be the characteristic function of the discrete sequence $u_j^n$ for a fixed $n=0,1,\dots.$ Then equation (\ref{fr_eq:time_space_frac_discr_eqn_0}), in terms of characteristic functions, takes the form
\begin{equation} 
  \hat{u}^{n+1} (\xi)=
  {\ds \gamma_n  + \sum_{m=2}^{n} c_m \hat{u}^{n+1-m}(\xi) +
                    \hat{q}(\xi) \hat{u}^n(\xi)}, 
  \label{char_050}
  \end{equation} 
  since $\hat{u}^{0}(\xi)=1.$ Further, let $\hat{U}_{\tau}(s,\xi)$ be the discrete Laplace transform of $\hat{u}^{n+1}(\xi),$ namely
  \[
  \hat{U}_{\tau}(s,\xi)=\tau \sum_{n=0}^{\infty} \hat{u}^{n+1}(\xi)e^{-s t_n}, \, s>0.
  \]
Then multiplying both sides of (\ref{char_050}) by $\tau e^{-n \tau s}$ and summing over the index $n,$ one obtains
  \[
\hat{U}_{\tau}(s, \xi) = \gamma_{\tau}(s) + \tau \sum_{n=0}^{\infty}  \left( \sum_{m=2}^{n+1} c_m \hat{u}^{n+1-m}(\xi) \right) e^{- n \tau s} + \hat{q}(\xi) \tau \sum_{n=0}^{\infty} \hat{u}^{n}(\xi) e^{-s n \tau}
\]
\[=
\gamma_{\tau}(s) - \tau \sum_{n=0}^{\infty}  \left( \sum_{m=1}^{n+1} (-1)^m {\beta \choose m} \hat{u}^{n+1-m}(\xi) \right) e^{- n \tau s} 
\]
\begin{equation}
+ \hat{d}(\xi) \tau^{1+\beta} \sum_{n=0}^{\infty} \hat{u}^{n}(\xi) e^{-s n \tau},
\label{char_051}
  \end{equation}
  where
  \[
  \gamma_{\tau}(s)= \tau \sum_{n=0}^{\infty} \gamma_n e^{-sn\tau}= \tau \sum_{n=0}^{\infty} \sum_{m=0}^{n+1} (-1)^m {\beta \choose m} e^{-sn\tau}.
  \]
  Changing the order of summation one can show that
   \[
  \gamma_{\tau}(s) = e^{s\tau} \left(\sum_{n=0}^{\infty} \tau e^{-sn\tau}\right)  \sum_{m=0}^{\infty} (-1)^m {\beta \choose m} e^{-sm\tau}.
  \] 
  In order to prove the theorem we need to show that $\hat{U}_{\tau}(s,h\xi)$ converges  as $h \to 0$ (that implies $\tau \to 0$ too) to 
 \[
 L[E_{\beta}(\Psi(\xi)t^{\beta})](s) = \frac{s^{\beta-1}}{s^\beta+\big(-\Psi(\xi)\big)},
 \] 
 the Laplace transform of the Mittag-Leffler function $E_{\beta}(x)$ composed by $\Psi(\xi)t^{\beta}.$ Indeed, this convergence implies the convergence  $\hat{u}^n(h\xi) \to E_{\beta}(\Psi(\xi)t^{\beta}),$ as $n \to \infty,$ uniformly for all $\xi \in \mathcal{K},$ where $\mathcal{K}$ is an arbitrary compact in $\re^d.$ In turn, the latter convergence is equivalent to the convergence in law of the sequence $S_n$ to the process $Y_{W_t}.$
To show the convergence $\hat{U}_{\tau}(s,h\xi) \to L[E_{\beta}(\Psi(\xi)t^{\beta})](s),$ we notice that
 \begin{equation}
 \label{char_052}
 \tau \sum_{n=0}^{\infty} \hat{u}^{n}(\xi) e^{-s n \tau} = \tau +e^{-s\tau} \hat{U}_{\tau}(s,\xi),
 \end{equation}
  and changing the order of summation
  \[
  \tau \sum_{n=0}^{\infty}   \left( \sum_{m=1}^{n+1} (-1)^m  {\beta \choose m} \hat{u}^{n+1-m}(\xi) \right) e^{- n \tau s}  
  \]
  \begin{equation}= -\tau \beta + \Big(\tau e^{s\tau}+\hat{U}_{\tau}(s,\xi)\Big) \left(\sum_{n=0}^{\infty} (-1)^n {\beta \choose n}e^{-sn\tau} -1 \right). 
  \label{char_053}
  \end{equation}
 It follows from equations (\ref{char_051})-(\ref{char_053}) that
  \begin{equation}
  \label{char_500}
  \hat{U}_{\tau}(s,\xi)= \frac{\ds{ e^{s\tau}  I_{\tau}(\beta,s) \left(\sum_{n=0}^{\infty} \tau e^{-sn\tau} - \tau \right) + \tau^{1-\beta} (\tau^{\beta} \hat{d}(\xi)+\beta+e^{s\tau})}}{\ds{ I_{\tau}(\beta,s) -\hat{d}(\xi) e^{-s\tau} }}
  \end{equation}
  where
  \begin{equation*}
  \label{char_550}
  I_{\tau}(\beta,s)=\frac{1}{\tau^{\beta}} \sum_{n=0}^{\infty} (-1)^n {\beta \choose n} e^{-sn\tau}.
    \end{equation*}
 Further, the following limits hold:
 \begin{equation}
 \label{lim_1}
\lim_{\tau \to 0} \frac{1}{\tau^{\beta}} \sum_{n=0}^{\infty} (-1)^n {\beta \choose n} {e^{-sn\tau}} =  \left(_{-\infty} \mathcal{D}^{\beta}_t e^{st} \right)_{|_{t=0}}=s^{\beta},
\end{equation}
\begin{equation}
\label{lim_2}
\lim_{\tau \to 0} \left(\sum_{n=0}^{\infty} \tau e^{-sn\tau} - \tau \right)  =s^{-1},
\end{equation}
\begin{equation}
\label{lim_3}
\lim_{h \to 0} \hat{d}(h \xi) =\Psi(\xi),
\end{equation}
\begin{equation}
\label{lim_4}
\lim_{\tau \to 0} \tau^{1-\beta} (\tau^{\beta} \hat{d}(\xi)+\beta+e^{s\tau}) =0,
\end{equation} 
  The relation (\ref{lim_1}) follows from the definition (\ref{GL_1}) with $f(t)=e^{st}, s>0,$ and the relation (\ref{ex}).   The relations (\ref{lim_2}) and (\ref{lim_4}) can be easily verified by direct calculation. The relation (\ref{lim_3})  is proved in Corollary to Lemma \ref{lem}.   
Now taking into account the relations (\ref{lim_1})-(\ref{lim_4}) it follows from (\ref{char_500}) that
  \[
  \lim_{h \to 0} \hat{U}_{\tau}(s,h\xi) = \frac{s^{\beta-1}}{s^{\beta}-\Psi(\xi)},
  \]
proving the theorem.

\section{Final remarks}

Theorem \ref{thm} extends to the case when the left hand side of equation
(\ref{ctrw_010}) is a time distributed fractional order
differential operator with a mixing measure $\mu$ whose support satisfies $ \supp \, \mu \subseteq [0,1]:$ 
\begin{equation}
\label{dode}
D_{\mu} u(t,x)=\int_0^1 D_{\ast}^{\beta} u(t,x) d \mu(\beta) = \Psi(D)u(t,x), \quad t>0, \, x \in \re^d,
\end{equation} 
where $\Psi(D)$ is a pseudo-differential operator with the symbol $\Psi(\xi)$ defined in (\ref{fr_psi_010}).
In this case for the left hand side of (\ref{dode}) we again have a discretization of the form (\ref{eqn_011}). Namely, we have
\begin{equation*}
  D_{\mu} u^n_j \approx 
    a(\tau)                                              
  \left( 
     u^{n+1}_j - c^{\ast}_1 u^n_j - \sum_{m=2}^n c^{\ast}_m u^{n+1-m}_j - \gamma^{\ast}_n u^0_j 
  \right),
  \label{eqn_0110} 
\end{equation*}
where
\begin{equation}
\label{dode_01}
a(\tau) =\int_0^1 a(\tau,\beta) d \mu (\beta), \, \, c_k^{\ast} =\int_0^1 c_k (\beta) d\mu(\beta), \, k=1,\dots, n, 
\end{equation}
\begin{equation}
\label{dode_02}
\gamma_n^{\ast} =\int_0^1 \gamma_n (\beta) d\mu(\beta), \, \, n=1,2,\dots.
\end{equation}
In equations (\ref{dode_01}) and (\ref{dode_02}) the  integrands $a(\tau, \beta), \, c_k(\beta)$ and $\gamma_n(\beta)$ are defined in (\ref{c_1}),(\ref{c_3}) or (\ref{c_4}),(\ref{c_6}) depending on whether the Gr\"unwald-Letnikov or quadrature approximation in paper \cite{Liu} is used for discretization of $D_{\ast}^{\beta}u(t,x)$ in (\ref{dode}). The detailed analysis of the corresponding CTRW including simulation models will be presented in a separate paper. 

We also note that condition (\ref{stabcond2}) takes the form
$$\tau \le \Big(\frac{2-2^{1-\beta}}{\Gamma(2-\beta)Q(h)}\Big)^{1 \over \beta}$$
if the non-Markovian probabilities are selected as in paper \cite{Liu,Andries}.
This condition as well as (\ref{stabcond2}) generalize
the well-known Lax's stability condition   $\tau \le h^2/2$  arising in the finite-difference method
for solution of an initial value problem for the heat equation, which corresponds to the case $\beta=1.$
In this case $Q(h)$ reduces simply to $Q(h)=2/h^2.$

Finally, in the particular case $\beta=1,$ Theorem \ref{thm} reduces to the following theorem, which provides a random walk approximation of stochastic
processes $X^{\rho}_t \in \mathbb{SS}.$ 

\begin{theorem} 
Let $X_j \in h\Z^d, \, j \ge 1,$ be i.i.d.
random vectors with
the probability mass function $p_{k} = \mathbb{P}(X_1=k)$ defined in (\ref{trpr111}) 
with some $\tau >0,$ $h>0.$ 
Assume that
\begin{equation*}
\label{fr_rw_052}
\sigma(\tau, h):= 2 \tau \sum_{m \neq 0}
\frac{Q_m(h)}{|m|^{n}} \leq 1.
\end{equation*}
Then the sequence of random vectors ${S}_{N}={X}_{1}+ ... +  {
X}_{N},$ converges  in law as
 $N \rightarrow \infty$ to $X^{\rho}_t \in \mathbb{SS},$ whose probability density function is $G^{\rho}(t,x)$ defined in (\ref{fr_sln_090}), that is the solution to
equation (\ref{equation}) with the initial condition $u(0,x)=\delta_0 (x).$ 
\end{theorem}

This theorem in the particular case $d\rho (\alpha) = a(\alpha)d\alpha,$ where $a(\cdot)$ is a positive continuous function on the interval $[0,2],$ is proved in \cite{UmarovSteinberg,Andries}. 


\end{document}